\documentclass{article}

\usepackage{latexsym,amssymb}

\makeatletter

 \twocolumn
\def\addcontentsline#1#2#3{}

\long\def\icmltitle#1{%
  {\center\baselineskip 18pt
                       \toptitlebar{\Large\bf #1}\bottomtitlebar}
}
\def\toptitlebar{\hrule height1pt \vskip .25in}
\def\bottomtitlebar{\vskip .22in \hrule height1pt \vskip .3in}
\def\icmlauthor#1#2{\par {\bf #1} \hfill {\sc #2}}
\long\def\icmladdress#1{\par\vskip 0.03in #1 \vskip 0.10in}

\renewenvironment{abstract}
   {%
\centerline{\large\bf Abstract}
    \vspace{-0.12in}\begin{quote}}
   {\par\end{quote}\vskip 0.12in}

\def\@startsection#1#2#3#4#5#6{\if@noskipsec \leavevmode \fi
   \par \@tempskipa #4\relax
   \@afterindenttrue
   \ifdim \@tempskipa <\z@ \@tempskipa -\@tempskipa \fi
   \if@nobreak \everypar{}\else
     \addpenalty{\@secpenalty}\addvspace{\@tempskipa}\fi \@ifstar
     {\@ssect{#3}{#4}{#5}{#6}}{\@dblarg{\@sict{#1}{#2}{#3}{#4}{#5}{#6}}}}

\def\@sict#1#2#3#4#5#6[#7]#8{\ifnum #2>\c@secnumdepth
     \def\@svsec{}\else
     \refstepcounter{#1}\edef\@svsec{\csname the#1\endcsname}\fi
     \@tempskipa #5\relax
      \ifdim \@tempskipa>\z@
        \begingroup #6\relax
          \@hangfrom{\hskip #3\relax\@svsec.~}{\interlinepenalty \@M #8\par}
        \endgroup
       \csname #1mark\endcsname{#7}\addcontentsline
         {toc}{#1}{\ifnum #2>\c@secnumdepth \else
                      \protect\numberline{\csname the#1\endcsname}\fi
                    #7}\else
        \def\@svsechd{#6\hskip #3\@svsec #8\csname #1mark\endcsname
                      {#7}\addcontentsline
                           {toc}{#1}{\ifnum #2>\c@secnumdepth \else
                             \protect\numberline{\csname the#1\endcsname}\fi
                       #7}}\fi
     \@xsect{#5}}

\def\@sect#1#2#3#4#5#6[#7]#8{\ifnum #2>\c@secnumdepth
     \def\@svsec{}\else
     \refstepcounter{#1}\edef\@svsec{\csname the#1\endcsname\hskip 0.4em }\fi
     \@tempskipa #5\relax
      \ifdim \@tempskipa>\z@
        \begingroup #6\relax
          \@hangfrom{\hskip #3\relax\@svsec}{\interlinepenalty \@M #8\par}
        \endgroup
       \csname #1mark\endcsname{#7}\addcontentsline
         {toc}{#1}{\ifnum #2>\c@secnumdepth \else
                      \protect\numberline{\csname the#1\endcsname}\fi
                    #7}\else
        \def\@svsechd{#6\hskip #3\@svsec #8\csname #1mark\endcsname
                      {#7}\addcontentsline
                           {toc}{#1}{\ifnum #2>\c@secnumdepth \else
                             \protect\numberline{\csname the#1\endcsname}\fi
                       #7}}\fi
     \@xsect{#5}}

\def\section{\@startsection{section}{1}{\z@}{-0.12in}{0.02in}
             {\large\bf\raggedright}}
\def\subsection{\@startsection{subsection}{2}{\z@}{-0.10in}{0.01in}
                {\normalsize\bf\raggedright}}
\def\subsubsection{\@startsection{subsubsection}{3}{\z@}{-0.08in}{0.01in}
                {\normalsize\sc\raggedright}}
\def\paragraph{\@startsection{paragraph}{4}{\z@}{1.5ex plus
  0.5ex minus .2ex}{-1em}{\normalsize\bf}}
\def\subparagraph{\@startsection{subparagraph}{5}{\z@}{1.5ex plus
  0.5ex minus .2ex}{-1em}{\normalsize\bf}}

\skip\footins 9pt
\def\footnoterule{\kern-3pt \hrule width 0.8in \kern 2.6pt }
\labelwidth\leftmargini\advance\labelwidth-\labelsep \labelsep 5pt

\belowdisplayskip \abovedisplayskip
\def\@normalsize{\@setsize\normalsize{11pt}\xpt\@xpt}
\def\small{\@setsize\small{10pt}\ixpt\@ixpt}
\def\footnotesize{\@setsize\footnotesize{10pt}\ixpt\@ixpt}
\def\scriptsize{\@setsize\scriptsize{8pt}\viipt\@viipt}
\def\tiny{\@setsize\tiny{7pt}\vipt\@vipt}
\def\large{\@setsize\large{14pt}\xiipt\@xiipt}
\def\Large{\@setsize\Large{16pt}\xivpt\@xivpt}
\def\LARGE{\@setsize\LARGE{20pt}\xviipt\@xviipt}
\def\huge{\@setsize\huge{23pt}\xxpt\@xxpt}
\def\Huge{\@setsize\Huge{28pt}\xxvpt\@xxvpt}

\newsavebox\captionbox\newdimen\captionboxwid

\long\def\@makecaption#1#2{
 \vskip 10pt
        \baselineskip 11pt
        \setbox\@tempboxa\hbox{#1. #2}
        \ifdim \wd\@tempboxa >\hsize
        \sbox{\captionbox}{\small\sl #1.~}
        \captionboxwid=\wd\captionbox
        \usebox\captionbox {\footnotesize #2}
        \else
          \centerline{{\small\sl #1.} {\small #2}}
        \fi}

\def\fnum@figure{Figure \thefigure}
\def\fnum@table{Table \thetable}

\def\texitem#1{\par\noindent\hangindent 12pt
               \hbox to 12pt {\hss #1 ~}\ignorespaces}

\long\def\comment#1{}

\makeatother

\def\npcite{\cite}
\newtheorem{Def1}{Definition}

\newtheorem{Theorem}[Def1]{Theorem}
\newtheorem{Cor}[Def1]{Corollary}

\newtheorem{Ex1}[Def1]{Example}

\newenvironment{Example} {\begin{Ex1} \begin{upshape}} {\end{upshape} \end{Ex1}}

\def\defaultskip{\medskip}

\def\beq{\begin{equation}}
\def\eeq{\end{equation}}
\def\beqn{\begin{displaymath}}
\def\eeqn{\end{displaymath}}
\def\bqa{\begin{eqnarray}}
\def\eqa{\end{eqnarray}}
\def\bqan{\begin{eqnarray*}}
\def\eqan{\end{eqnarray*}}

\def\calC{\mathcal C}

\def\calF{\mathcal F}

\def\calX{\mathcal X}
\def\calY{\mathcal Y}

\def\NNN{\mathbb N}
\def\RRR{\mathbb R}
\def\QQQ{\mathbb Q}

\def\Expect{{\mathbf E}}
\def\Prob{{\mathbf P}}
\def\rrho{\varrho}
\def\eps{\varepsilon}

\def\leqt{_{1:t}}
\def\ltt{_{<t}}
\def\leqn{_{1:n}}
\def\ltn{_{<n}}

\def\ltinf{_{<\infty}}
\def\_norm{_\mathrm{norm}}

\def\for_all{\mbox{ for all }}
\def\such_that{\mbox{ such that }}

\def\und{\mbox{ and }}

\def\lb{{\log_2}}                      % binary logarithm
\def\th{\vartheta}

\def\toinfty#1{\stackrel{#1\to\infty}{\longrightarrow}}

\def\subsection#1{\vspace{1ex}\noindent{\bf{#1.}}}

\def\_norm{_\mathrm{n}}
\def\Clin{\calC^\mathrm{lin1}}
\def\Fgauss{\calF^\mathrm{Gauss}}

\def\xnltn{x_{1:n}}
\def\xtltt{x_{1:t}}

\def\citet{\npcite}
\def\citep{\cite}
\def\citeauthor{\aucite}

\begin{document}

\twocolumn[IDSIA-02-05 \hfill 16 January 2005
\icmltitle{Strong Asymptotic Assertions for Discrete MDL \\ in Regression and Classification}
\icmlauthor{Jan Poland}{jan@idsia.ch}
\icmlauthor{Marcus Hutter}{marcus@idsia.ch}
\icmladdress{IDSIA, Galleria 2, CH-6928 Manno-Lugano, Switzerland \hfill \sc www.idsia.ch}
\vskip 0.3in
]

\begin{abstract}
We study the properties of the MDL (or maximum penalized
complexity) estimator for Regression and Classification, where the
underlying model class is countable. We show in particular a
finite bound on the Hellinger losses under the only assumption
that there is a ``true" model contained in the class. This implies
almost sure convergence of the predictive distribution to the true
one at a fast rate. It corresponds to Solomonoff's central theorem
of universal induction, however with a bound that is exponentially
larger.

\end{abstract}

%%%%%%%%%%%%%%%%%%%%%%%%%%%%%%%%%%%%%%%%%%%%%%%%%%%%%%%%%%%%%%%
\section{Introduction}
%%%%%%%%%%%%%%%%%%%%%%%%%%%%%%%%%%%%%%%%%%%%%%%%%%%%%%%%%%%%%%%

Bayesian methods are popular in Machine Learning. So it is natural
to study their predictive properties: How do they behave
asymptotically for increasing sample size? Are loss bounds
obtainable, either for certain interesting loss functions or even
for more general classes of loss functions?

In this paper, we consider the two maybe most important Bayesian
methods for prediction in the context of regression and
classification. The first one is \emph{marginalization}: Given
some data and a model class, obtain a predictive model by
integrating over the model class. This \emph{Bayes mixture} is
``ideal" Bayesian prediction in many respects, however in many
cases it is computationally untractable. Therefore, a commonly
employed method is to compute a \emph{maximum penalized
complexity} or \emph{maximum a posteriori (MAP)} or \emph{minimum
description length (MDL)} estimator. This predicts according to
the ``best" model instead of a mixture. The MDL principle is
important for its own sake, not only as approximation of the Bayes
mixture.

Most work on Bayesian prediction has been carried out for
\emph{continuous} model classes, e.g.\ classes with one free
parameter $\th\in\RRR^d$. While the predictive properties of the
Bayes mixture are excellent under mild conditions
\cite{Clarke:90,Hutter:03optisp,Ghosal:00,Hutter:04uaibook}, corresponding MAP or
MDL results are more difficult to establish. For MDL in the strong
sense of
\emph{description length}, the parameter space has to be
discretized appropriately (and dynamically with increasing sample
size) \cite{Rissanen:96,Barron:98,Barron:91}. A MAP estimator on
the other hand can be very bad in general. In statistical
literature, some important work has been performed on the
asymptotical discovery of the true parameter, e.g.\
\cite{LeCam:00}. This can only hold if each model occurs no more
than once in the class. Thus it is violated e.g.\ in the case of
an artificial neural network, where exchanging two hidden units in
the same layer does not alter the network behavior.

In the case of \emph{discrete} model classes, both loss bounds
and asymptotic assertions for the Bayes mixture are relatively
easy to prove, compare Theorem \ref{thSolomonoff}. In
\cite{Poland:04mdl}, corresponding results for MDL were shown.
The setting is sequence prediction but otherwise very general.
The only assumption necessary is that the true distribution is
contained in the model class. Assertions are given directly
for the predictions, thus there is no problem of possibly
undistinguishable models. In order to prove that the MDL
estimator (precisely, the \emph{static} MDL estimator in terms
of \cite{Poland:04mdl}) has good predictive properties, we
introduce an intermediate step and show first the predictive
properties of \emph{dynamic MDL}, where a new MDL estimator is
computed for each possible next observation.

In this paper, we will derive analogous results for regression and
classification. While results for classification can be generalized
from sequence prediction by conditionalizing everything to the
input, regression is technically more difficult.
Therefore the next section, which deals with the regression setup,
covers the major part of the paper. Instead of the popular
Euclidian and Kullback-Leibler distances for measuring prediction
quality we need to exploit the Hellinger distance. We show that
online MDL converges to the true distribution in mean Hellinger
sum, which implies ``rapid'' convergence with probability one.
Classification is briefly discussed in Section \ref{secClass},
followed by a discussion and conclusions in Section \ref{secConc}.

%%%%%%%%%%%%%%%%%%%%%%%%%%%%%%%%%%%%%%%%%%%%%%%%%%%%%%%%%%%%%%%
\section{Regression}\label{secReg}
%%%%%%%%%%%%%%%%%%%%%%%%%%%%%%%%%%%%%%%%%%%%%%%%%%%%%%%%%%%%%%%

We neglect computational aspects and study the properties of the
\emph{optimal} Bayes mixture and MDL predictors.
When a new sample is observed, the estimator is updated. Thus,
regression is considered in an \emph{online framework}: The
first input $x_1$ is presented, we predict the output $y_1$
and then observe its true value, the second input $x_2$ is
presented and so on.

%--------------------------------------------------------------
\subsection{Setup}\label{ssecSetup}
%--------------------------------------------------------------
Consider a regression problem with arbitrary domain $\calX$ (we
need no structural assumptions at all on $\calX$) and co-domain
$\calY=\RRR$. The task is to learn/fit/infer a function
$f:\calX\to\calY$, or more generally a conditional probability
density $\nu(y|x)$, from data $\{(x_1,y_1),...,(x_n,y_n)\}$.
Formally, we are given a countable class $\calC$ of models that
are functions $\nu$ from $\calX$ to \emph{uniformly bounded
probability densities} on $\RRR$. That is, $\calC=\{\nu_i:\,i\geq
1\}$, and there is some $C>0$ such that
\bqa
\label{eqC}
\lefteqn{ 0\leq\nu_i(y|x)\leq C \und
\int_{-\infty}^\infty \nu_i(y|x)dy=1}\\ \nonumber
&&\for_all i\geq 1,\ x\in\calX, \und y\in\calY.
\eqa
Each $\nu$ induces a probability density on
$\RRR^n$ for $n$-tuples $x\leqn\in\calX^n$ by
$\nu(y\leqn|x\leqn)=\prod_{t=1}^n\nu(y_t|x_t)$. The notation
$x\leqn$ for $n$-tuples is common in sequence prediction. Each
model $\nu\in\calC$ is associated with a \emph{prior weight}
$w_\nu>0$. The logarithm $\lb w_\nu$ has often an interpretation
as model \emph{complexity}. We require $\sum_\nu w_\nu=1$. Then by
the Kraft inequality, one can assign to each model $\nu\in\calC$ a
prefix-code of length $^\lceil\lb w_\nu\!^\rceil$.

We assume that an infinite stream of data
$(x_{1:\infty},y_{1:\infty})$ % $(x\ltinf,y\ltinf)$ undefined
is generated as follows: Each $x_t$ may be produced by an
arbitrary mechanism, while $y_t$ is sampled from a \emph{true
distribution} $\mu$ conditioned on $x_t$. In order to obtain strong
convergence results, we will require that $\mu\in\calC$.

\begin{Example} \label{ex1} Take $\calX=\RRR$ and
$\Clin_\sigma\cong\{ax+b+N(0,\sigma^2):a,b\in\QQQ\}$ to be the class of linear
regression models with rational coefficients $a$, $b$, and independent
Gaussian noise of fixed variance $\sigma^2>0$. That is,
$\Clin_\sigma=\{\nu^{a,b,\sigma}:a,b\in\QQQ\}$, where
\beqn
\nu^{a,b,\sigma}(x,y)=\phi_{\sigma^2}(y-ax-b)={\textstyle\frac{1}{\sqrt{2\pi\sigma^2}}}
e^{-\frac{1}{2\sigma^2}(y-ax-b)^2}.
\eeqn
Alternatively, you may consider the class
$\Clin_{\geq\sigma_0}=\{\nu^{a,b,\sigma}:a,b,\sigma\in\QQQ,\sigma\geq\sigma_0\}$
for some $\sigma_0>0$, where also the noise amplitude is part of
the models. In the following, we also discuss how to admit
degenerate Gaussians that are point measures such as
$\Clin_{\geq 0}$.
\end{Example}

The setup (\ref{eqC}) guarantees that all subsequent MDL
estimators [(\ref{eqStaticMDL}) and (\ref{eqDynamicMDL})] exist.
However, our results and proofs generalize in several directions.
First, for the co-domain $\calY$ we may choose any $\sigma$-finite
measure space instead of $\RRR$, since we need only Radon-Nikodym
densities below. Second, the uniformly boundedness condition can
be relaxed, if the MDL estimators still exist. This holds for
example for the class $\Clin_{\geq 0}$ (see the preceding
example), if the definition of the MDL estimators is adapted
appropriately (see footnote \ref{footnoteDirac} on page
\pageref{footnoteDirac}). Third, the results remain valid for
semimeasures with $\int\nu\leq 1$ instead of measures and $\sum
w_\nu\leq 1$, which is however not very relevant for regression
(but for universal sequence prediction). In order to keep things
simple, we develop all results on the basis of (\ref{eqC}). Note
finally that the models in $\calC$ may be time-dependent, and we
need not even make this explicit, since the time can be
incorporated into $\calX$
($x_t=(x'_t,t)\in\calX'\times\NNN=\calX$). In this way we may also
make the models depend on the actual past outcome, if this is
desired ($x_t=(x'_{1:t}y_{1:t-1})\in\calX'^*\times\calY^*=\calX$).

The case of independent Gaussian noise as in Example \ref{ex1} is
a particularly important one. We therefore introduce the family
\bqa
\label{eqFgauss}
\Fgauss_\geq &=&
\Big\{\calC=\{\nu_i,\sigma_i\}_{i=1}^\infty:\nu_i(x,y)=
\\ \nonumber&&\phi_{\sigma_i^2}\big(y-f_i(x)\big),
\sigma_i\geq\sigma_0>0,f_i\!:\!\calX\!\to\!\RRR\Big\}.
\eqa
of all countable regression model classes with lower bounded
Gaussian noise. Clearly,
$\Clin_\sigma,\Clin_{\geq\sigma_0}\in\Fgauss_\geq$ is satisfied.
Similarly $\Fgauss\supset\Fgauss_\geq$ denotes the corresponding
family without lower bound on $\sigma_i$. Then
$\Clin_{\geq 0}\in\Fgauss\setminus\Fgauss_\geq$.

We define the \emph{Bayes mixture}, which for each
$n\geq 1$ maps an $n$-tuple of inputs
$x\leqn\in\calX^n$ to a probability density on $\RRR^n$:
\beq
\label{eqxi}
\xi(y\leqn|x\leqn) = \sum_{\nu\in\calC}w_\nu\nu(y\leqn|x\leqn)=
\sum_{\nu\in\calC}w_\nu\prod_{t=1}^{n}\nu(y_t|x_t)
\eeq
(recall $\sum_\nu w_\nu=1$). Hence, the Bayes mixture
\emph{dominates} each $\nu$ by means of
$\xi(\cdot|x\leqn)\geq w_\nu \nu(\cdot|x\leqn)$ for all $x\leqn$.
For $\nu\in\calC$ and $x_n\in\calX$, the $\nu$-prediction of
$y_n\in\RRR$, that is the $nu$-probability density of observing $y_n$,
is
$$\nu(y_n|\xnltn,y\ltn)=\nu(y_n|x_n).$$
This is independent of the history
$(x\ltn,y\ltn)=(x_{1:n-1},y_{1:n-1})$. In contrast, the \emph{Bayes
mixture prediction} or regression, which is also a measure on
$\RRR$, depends on the history:
\iftrue
\beq
\label{eqxipred}
\xi(y_n|\xnltn,y\ltn)=\frac{\xi(y\leqn|x\leqn)}{\xi(y\ltn|x\ltn)}
=\frac{\sum_\nu w_\nu\prod_{t=1}^{n}\nu(y_t|x_t)} {\sum_\nu
w_\nu\prod_{t=1}^{n-1}\nu(y_t|x_t)}.
\eeq
\else
\beq
\label{eqxipred}
\xi(y_n|\xnltn,y\ltn)=\frac{\xi(y\leqn|x\leqn)}{\xi(y\ltn|x\ltn)}
=\sum_\nu w_\nu^n\prod_{t=1}^{n}\nu(y_t|x_t), \quad\mbox{where}\quad
w_\nu^n:=\frac{w_\nu} {\sum_\nu w_\nu\prod_{t=1}^{n-1}\nu(y_t|x_t)}
\eeq
are updated weights depending on past data $(x\ltn,y\ltn)$.
\fi
This is also known as \emph{marginalization}. Observe that the
denominator in (\ref{eqxipred}) vanishes only on a set of
$\mu$-measure zero, if the true distribution
$\mu$ is contained in $\calC$. Under condition (\ref{eqC}), the
Bayes mixture prediction is uniformly bounded. It can be argued
intuitively that in case of unknown
$\mu\in\calC$ the Bayes mixture is the best possible model for
$\mu$. Formally, its predictive properties are excellent:

\begin{Theorem}
\label{thSolomonoff}
Let $\mu\in\calC$, $n\geq 1$, and $x\leqn\in\calX^n$, then
\bqa
\label{eqSolomonoff}
\sum_{t=1}^n \Expect \int
  \Big( \sqrt{\mu(y_t|\xtltt,y\ltt)}-\sqrt{\xi(y_t|\xtltt,y\ltt)} \Big)^2 dy_t \\ \nonumber \ \leq \ \ln w_\mu^{-1}.
\eqa
\end{Theorem}

$\Expect$ denotes the expectation with respect to the true distribution
$\mu$. Hence in this case we have $\Expect\ldots=\int\ldots
\mu(dy\ltt)$. The integral expression is also known as
\emph{square Hellinger distance}.
It will emerge as a main tool in the subsequent proofs. So the
theorem states that on any input sequence $x\ltinf$ the
expected cumulated Hellinger divergence of
$\mu$ and the Bayes mixture prediction is bounded by $\ln w_\mu^{-1}$.
A closely related result was discovered by Solomonoff
(\cite{Solomonoff:78}) for universal sequence prediction, a
``modern" proof can be found in \cite{Hutter:04uaibook}. This
proof can be adapted in our regression framework. Alternatively,
it is not difficult to give a proof in a few lines analogous to
(\ref{eqKLbound1}) and (\ref{eqKLbound2}) by using
(\ref{eqHleqKL}).

We introduce the term \emph{convergence in mean Hellinger sum (i.m.H.s.)}
for bounds like (\ref{eqSolomonoff}): For some predictive density
$\psi$, the $\psi$-predictions converge to the $\mu$-predictions
i.m.H.s.\
%(short $\psi\toiHs\mu$)
on a sequence of inputs $x\ltinf\in\calX^\infty$, if there is
$R>0$ such that $H^2_{x\ltinf}(\mu,\psi)\leq R$, where
\bqa
\label{eqHellinger}
\lefteqn{H^2_{x\ltinf}(\mu,\psi) = \sum_{t=1}^\infty \Expect [h^2_t]
\mbox{ with }}\\ \nonumber &&h^2_t=\int \Big(
\sqrt{\mu(y_t|\xtltt,y\ltt)}-\sqrt{\psi(y_t|\xtltt,y\ltt)}
\Big)^2dy_t.
\eqa
Convergence i.m.H.s.\ is a very strong convergence criterion. It
asserts a finite expected cumulative Hellinger loss in the first
place. If the co-domain $\calY$ is finite as for classification
(see Section \ref{secClass}), then convergence i.m.H.s.\ implies
almost sure (a.s.) convergence of the (finitely many) posterior
probabilities. For regression, the situation is more complex,
since the posterior probabilities are densities, i.e.\ Banach space
valued. Here, convergence i.m.H.s.\ implies that with
$\mu$-probability one the square roots of the predictive densities
converge to the square roots of the
$\mu$-densities in $L^2(\RRR)$ (endowed with the Lebesgue
measure). In other words, $h^2_t$ converges to zero a.s.:
\bqa
\label{eqIHSWP1}
\!\!\!\Prob\Big(\exists t\geq n: h^2_t\geq\eps\Big)&=&
\Prob\Big(\bigcup_{t\geq
n}\big\{h^2_t\geq\eps\big\}\Big)\\\nonumber&\leq&
\sum_{t\geq n}\Prob \big(h^2_t\geq\eps\big)
\\\nonumber&\leq&\frac{1}{\eps}\sum_{t=n}^\infty\Expect h^2_t\toinfty n 0
\eqa
holds by the union bound, the Markov inequality for all $\eps>0$,
and $H^2_{x\ltinf}<\infty$, respectively, where $\Prob$ is the
$\mu$-probability. If the densities are uniformly bounded, then
also the differences of the densities (as opposed to the
difference of the square roots) converge to zero:
\beqn
  \psi(y_t|\xtltt,y\ltt)-\mu(y_t|\xtltt,y\ltt)
  \stackrel{t\to\infty}\longrightarrow 0
  \quad\mbox{in $L^2(\RRR)$ a.s.}
\eeqn
Moreover, the finite bound on the cumulative
Hellinger distances can be interpreted as a convergence rate.
Compare the parallel concept ``convergence in mean sum"
\cite{Hutter:03optisp,Poland:04mdl,Hutter:04uaibook}.

%--------------------------------------------------------------
\subsection{MDL Predictions}\label{ssecMDL}
%--------------------------------------------------------------
In many cases, the Bayes mixture is not only intractable, but even
hard to approximate. So a very common substitute is the (ideal)
MDL%
\footnote{\label{footnoteMDL} There is some disagreement about the
exact meaning of the term MDL. Sometimes a specific prior is
associated with MDL, while we admit arbitrary priors. More
importantly, when coding some data
$x$, one can exploit the fact that once the model $\nu^*$ is
specified, only data which lead to the maximizing element
$\nu^*$ need to be considered. This allows for a shorter
description than $\lb
\nu^*(x)$. Nevertheless, the \emph{construction principle} is
commonly termed MDL, compare for instance the ``ideal MDL" in
\cite{Vitanyi:00}.} estimator, also known as maximum a
posteriori (MAP) or maximum complexity penalized likelihood
estimator. Given a model class
$\calC$ with weights $(w_\nu)$ and a data set $(x\leqn,y\leqn)$,
we define the two-part MDL estimator as
\bqa
\nonumber
\nu^* \ = \ \nu^*_{(x\leqn,y\leqn)} & = &\arg\max_{\nu\in\calC}\{w_\nu
\nu(y\leqn|x\leqn)\}\und\\
\label{eqMDL}
  \rrho(y\leqn|x\leqn) & = & \max_{\nu\in\calC}\{w_\nu \nu(y\leqn|x\leqn)\}
  \\ \nonumber &=&w_{\nu^*}\nu^*(y\leqn|x\leqn).
\eqa
Note that we define both the model $\nu^*$ which is the MDL
estimator and its weighted density $\rho$. In our setup
(\ref{eqC}), the MDL estimator is well defined, since all maxima exist%
\footnote{\label{footnoteDirac}%
For a model class with Gaussian noise $\calC\in\Fgauss$
(\ref{eqFgauss}), we may dispose of the uniform boundedness
condition and admit e.g. also $\Clin_{\geq 0}$. In order to
compute the MDL estimator, we must then first check if there is
nonzero mass concentrated on $(x\leqn,y\leqn)$, in which case the
mass is even one and the corresponding model with the largest
weight is chosen. Otherwise, the MDL estimator is chosen according
to the maximum penalized density. All results and proofs below
generalize to this case.}. Moreover,
$\rrho(\cdot|x\leqn)$ is a density but its integral is less than 1
in general. We have
$\rrho(\cdot|x\leqn)\geq w_\nu \nu(\cdot|x\leqn)$, so like $\xi$,
$\rrho$ dominates each $\nu\in\calC$. Also,
$\rrho(\cdot|x\leqn)\leq \xi(\cdot|x\leqn)$ is clear by
definition. If we use $\nu^*$ for (sequential online)
prediction, this is the \emph{static MDL prediction}:
\beq
\label{eqStaticMDL}
  \rrho^{\mathrm{static}}(y_n|\xnltn,y\ltn)
  = \nu^*_{(x\ltn,y\ltn)}(y_n|x_n).
\eeq
This is the common way of using MDL for prediction. Clearly, the
static MDL predictor is a probability density on $\RRR$.
Alternatively, we may compute the MDL estimator for each possible
$y_n$ separately, arriving at the \emph{dynamic} MDL predictor:
\beq
\label{eqDynamicMDL}
  \rrho(y_n|\xnltn,y\ltn) = \frac{\rrho(y\leqn|x\leqn)}{\rrho(y\ltn|x\ltn)}.
\eeq
We have $\rrho(y_n|\xnltn,y\ltn) %$\frac{\rrho(y\leqn|x\leqn)}{\rrho(y\ltn|x\ltn)}
\leq \nu^*_{(x\leqn,y\leqn)}(y_n|x_n)$ for each $y_n$, which shows
that under condition (\ref{eqC}) the dynamic MDL predictor is
uniformly bounded. On the other hand,
$\rrho(y_n|\xnltn,y\ltn) %\frac{\rrho(y\leqn|x\leqn)}{\rrho(y\ltn|x\ltn)}
\geq \nu^*_{(x\ltn,y\ltn)}(y_n|x_n)$ holds, so the dynamic MDL
predictor may be a density with mass more than 1. Hence we must
usually \emph{normalize} it for predicting:
\beq
\label{eqNormDynamicMDL}
  \bar\rrho(y_n|\xnltn,y\ltn) = \frac{\rrho(y\leqn|x\leqn)}{\int\rrho(y\leqn|x\leqn)dy_n}.
\eeq
Both fractions in (\ref{eqDynamicMDL}) and
(\ref{eqNormDynamicMDL}) are well-defined except for a set of
measure zero. Dynamic MDL predictions are in a sense
computationally (almost) as expensive as the full Bayes mixture.

%--------------------------------------------------------------
\subsection{Convergence Results}\label{secConv}
%--------------------------------------------------------------
Our principal aim is to prove predictive properties of
\emph{static} MDL, since this is the practically most relevant
variant. To this end, we first need to establish corresponding
results for the dynamic MDL. Precisely, the following holds.

\begin{Theorem}
\label{th1}
Assume the setup (\ref{eqC}). If
$\mu\in\calC$, where $\mu$ is the true distribution, and $H^2_{x\ltinf}(\cdot,\cdot)$
is defined as in (\ref{eqHellinger}), then for all input sequences
$_{x\ltinf}\in\calX^\infty$ we have
\bqan
(i) && H^2_{x\ltinf}(\mu,\bar\rrho)\leq w_\mu^{-1}+\ln w_\mu^{-1}, \\
(ii) && H^2_{x\ltinf}(\bar\rrho,\rrho)\leq 2w_\mu^{-1}, \und \\
(iii) && H^2_{x\ltinf}(\rrho,\rrho^\mathrm{static})\leq
3w_\mu^{-1}.
\eqan
\end{Theorem}

Since the triangle inequality holds for
$\sqrt{H^2_{x\ltinf}(\cdot,\cdot)}$, we immediately conclude:

\begin{Cor}
\label{cor1} Given the setup (\ref{eqC}) and $\mu\in\calC$, then
all three predictors $\bar\rrho$, $\rrho$, and
$\rrho^\mathrm{static}$ converge to the true density $\mu$ in mean
Hellinger sum, for any input sequence $x\ltinf$. In particular, we
have $H^2(\mu,\rrho^\mathrm{static})\leq 21w_\mu^{-1}$.
\end{Cor}

We will only prove $(i)$ of Theorem \ref{th1} here. The proofs of
$(ii)$ and $(iii)$ can be similarly adapted from
\cite[Theorems 10 and 11]{Poland:04mdl}, since the Hellinger
distance is bounded by the absolute distance:
$\int\big( \sqrt{\mu(y)}-\sqrt{\nu(y)} \big)^2 dy\leq
\int \big|\mu(y)-\nu(y)\big|dy$ follows from
$(\sqrt a-\sqrt b)^2\leq|a-b|$ for any $a,b\in\RRR$ (this shows also
that the integral $h_t^2$ in (\ref{eqHellinger}) exists). In order
to show $(i)$, we make use of the fact that the squared Hellinger
distance is bounded by the
\emph{Kullback-Leibler divergence}:
\beq
\label{eqHleqKL}
\int
  \Big( \sqrt{\mu(y)}-\sqrt{\nu(y)} \Big)^2 dy\leq
\int \mu(y)\ln \frac{\mu(y)}{\nu(y)}dy
\eeq
for any two probability densities $\mu$ and $\nu$ on $\RRR$ (see
e.g. \cite[p. 178]{Borovkov:98}). So we only need to establish the
corresponding bound for the Kullback-Leibler divergence and show
\bqa
\nonumber
D_x(\mu\|\bar\rrho)&\!\!\!\!:=&\!\!\!\!
\sum_{t=1}^n \Expect
\int \mu(y_t|\xtltt,y\ltt)\ln
\frac{\mu(y_t|\xtltt,y\ltt)}{\bar\rrho(y_t|\xtltt,y\ltt)}dy_t\\
\label{eqKLbound}
&\leq& w_\mu^{-1}+\ln w_\mu^{-1}
\eqa
for all $n\geq 1$. In the following computation, we take $x\ltinf$
to be fixed and suppress it in the notation, writing e.g.
$\mu(y_t|y\ltt)$ instead of $\mu(y_t|\xtltt,y\ltt)$. Then
\bqa
\label{eqKLbound1}
D_x(\mu\|\bar\rrho)&=&
\sum_t \Expect\ln
\frac{\mu(y_t|y\ltt)}{\bar\rrho(y_t|y\ltt)}\\\nonumber&=&
\sum_t\Expect \left[
  \ln\frac{\mu(y_t|y\ltt)}{\rrho(y_t|y\ltt)}
  + \ln
  \frac{\int\rrho(y\leqt)dy_t}{\rrho(y\ltt)}\right].
\eqa
The first part of the last term is bounded by
\bqa
\label{eqKLbound2}
  \sum_t\Expect \ln\frac{\mu(y_t|y\ltt)}{\rrho(y_t|y\ltt)}
  &=& \Expect\ \ln\prod_{t=1}^n \frac{\mu(y_t|y\ltt)}{\rrho(y_t|y\ltt)}
  \\\nonumber&=& \Expect\ \ln\frac{\mu(y\leqn|x\leqn)}{\rrho(y\leqn|x\leqn)}
  \\\nonumber&\leq& \ln w_\mu^{-1},
\eqa
since always $\frac{\mu}{\rrho}\leq w_\mu^{-1}$. For the second
part, use $\ln u\leq u-1$ to obtain
\bqan
\lefteqn{\Expect \ln
\frac{\int\rrho(y\leqt)dy_t}{\rrho(y\ltt)}}\\
& \leq &
\!\!\!\!\sum_t \Expect
\left[\frac{\int\rrho(y\leqt)dy_t}{\rrho(y\ltt)}-1\right]\\
& =& \int \frac{\mu(y\ltt)(\int\rrho(y\leqt)dy_t-
\rrho(y\ltt)}{\rrho(y\ltt)}dy\ltt\\
&\leq&\!\!\!\! w_\mu^{-1}\left[\int\rrho(y\leqt)dy\leqt-
\int\rrho(y\ltt)dy\ltt\right].
\eqan
If this is summed over $t=1\ldots n$, the last term is
telescoping. So using $\rrho(\emptyset)=\max_\nu w_\nu\geq 0$ and
$\rrho\leq\xi$, we conclude
\bqa
\nonumber
\sum_t\Expect \ln \frac{\int\rrho(y\leqt)dy_t}{\rrho(y\ltt)}
&\leq &w_\mu^{-1}\left[\int\rrho(y\leqn)dy\leqn-
\rrho(\emptyset)\right]
\\
\label{eqKLbound3}&\leq& w_\mu^{-1}\int\xi(y\leqn)dy\leqn\\
\nonumber&=&w_\mu^{-1}.
\eqa
Hence, (\ref{eqKLbound1}), (\ref{eqKLbound2}), and
(\ref{eqKLbound3}) show together (\ref{eqKLbound}).
\hfill $\Box$

We may for example apply the result for the static predictions in
a Gaussian noise class $\calC\in\Fgauss$.

\begin{Cor} \label{cor2}
Let $\calC\in\Fgauss_\geq$ [see (\ref{eqFgauss})] then the mean
and the variance of the static MDL predictions converge to their
true values almost surely. The same holds for $\calC\in\Fgauss$.
In particular, if the variance of all models in $\calC$ is the
same value $\sigma^2$, then $\sum_t
2\big[1-\exp(-\frac{(g^*(x_t|\ldots)-f(x_t))^2}{8\sigma^2})\big]
\leq 21 w_\mu^{-1}$, where $f(x_t)$ is the mean value of the true
distribution and $g^*=\arg\min_{f_i}\{{1\over
n-1}\sum_{t=1}^{n-1}(y_t\!-\!f_i(x_t))^2+2\sigma^2\ln w_i^{-1}\}$ is
the mean of the MDL predictor.
\end{Cor}

For $\calC\in\Fgauss_\geq$, almost sure convergence holds since
otherwise the cumulative Hellinger distances would be infinite,
see (\ref{eqIHSWP1}). This generalizes to $\calC\in\Fgauss$;
compare the footnote \ref{footnoteDirac} on page
\pageref{footnoteDirac}. In the case of constant variance, the
cumulative Hellinger distances can be explicitly stated as above.
Note that since
$1-\exp\mbox{$\big(-\frac{(g^*(x_t|\ldots)-f(x_t))^2}{8\sigma^2}\big)$}\approx
\frac{(g^*(x_t|\ldots)-f(x_t))^2}{8\sigma^2}$ for small
$(g^*(x_t|\ldots)-f(x_t))^2$, this implies convergence of $g^*$ to
$f$ faster than $O(\frac{1}{\sqrt t})$ if the convergence is
monotone. Moreover, deviations of a fixed magnitude can only occur
finitely often.

Compared with the bound for the Bayes mixture in Theorem
\ref{thSolomonoff}, MDL bounds are exponentially larger. The
bounds are sharp, as shown in \cite[Example 9]{Poland:04mdl}, this
example may be also adapted to the regression framework.

%%%%%%%%%%%%%%%%%%%%%%%%%%%%%%%%%%%%%%%%%%%%%%%%%%%%%%%%%%%%%%%
\section{Classification}\label{secClass}
%%%%%%%%%%%%%%%%%%%%%%%%%%%%%%%%%%%%%%%%%%%%%%%%%%%%%%%%%%%%%%%

The classification setup is technically easier, since only a
finite co-domain $\calY$ has to be considered. Results
corresponding to Theorem \ref{th1} and Corollary \ref{cor1} follow
analogously. Alternatively, one may conditionalize the results for
sequence prediction in \cite{Poland:04mdl} with
respect to the input sequence $x\ltinf$, arriving equally at the
assertions for classification. The results in \cite{Poland:04mdl}
are formulated in terms of mean (square) sum convergence instead
of Hellinger sum convergence. On finite co-domain, these two
convergence notions induce the same topology.

\begin{Theorem}
Let $\calX$ be arbitrary and $\calY$ be a finite set of class
labels. $\calC=\{\nu_i:i\geq 1\}$ consists of classification
models, i.e.\ for each $\nu\in\calC$, $x\in\calX$ and
$y\in\calY$ we have $\nu(y|x)\geq 0$ and $\sum_y
\nu(y|x)=1$. Each model $\nu$ is associated with a prior weight $w_\nu>0$,
and $\sum_\nu w_\nu=1$ holds. Let the MDL predictions be defined
analogously to (\ref{eqMDL}), (\ref{eqStaticMDL}) and
(\ref{eqDynamicMDL}) (the difference being that here probabilities
are maximized instead of densities). Assume that $\mu\in\calC$,
where
$\mu$ is the true distribution. Then for each
$x\ltinf\in\calX^\infty$,
\bqan
\sum_{t=1}^\infty \Expect \sum_{y\in\calY}
\Big(\sqrt{\mu(y|x_t)}-\sqrt{\rrho^\mathrm{static}(y|\xtltt,y\ltt)}\Big)^2
&\!\!\!\!\!\!\!\! \leq &\!\!\!\! 21 w_\mu^{-1},\\
\sum_{t=1}^\infty \Expect \sum_{y\in\calY}
\Big(\mu(y|x_t)-\rrho^\mathrm{static}(y|\xtltt,y\ltt)\Big)^2
&\!\!\!\!\!\!\!\! \leq &\!\!\!\! 21 w_\mu^{-1}
\eqan
holds. Similar assertions are satisfied for the normalized and the
un-normalized dynamic MDL predictor. In particular, the predictive
probabilities of all three MDL predictors converge to the true
probabilities almost surely.
\end{Theorem}

The second bound on the quadratic differences is shown in
\cite{Poland:04mdl}. The assertions about almost sure convergence
follows as in (\ref{eqIHSWP1}).

%%%%%%%%%%%%%%%%%%%%%%%%%%%%%%%%%%%%%%%%%%%%%%%%%%%%%%%%%%%%%%%
\section{Discussion and Conclusions}\label{secConc}
%%%%%%%%%%%%%%%%%%%%%%%%%%%%%%%%%%%%%%%%%%%%%%%%%%%%%%%%%%%%%%%

We have seen that discrete MDL has good asymptotic predictive
properties. On the other hand, the loss bounds for MDL are
exponential compared to the Bayes mixture loss bound. This is
no proof artifact, as examples are easily constructed where
the bound is sharp \cite{Poland:04mdl}.

This has an important implication for the practical use of
MDL: One need to choose the underlying model class and the
prior carefully. Then it can be expected that the predictions
are good and converge fast: this is supported by theoretical
arguments in \cite{Rissanen:96,Poland:04mdlspeed}. The Bayes
mixture in contrast, which can be viewed as a very large
(infinite) weighted committee, also converges rapidly with
unfavorable model classes, but at higher computational
expenses.

One might be interested in other loss functions than the
Hellinger loss. For the classification case, a bound on the
expected error loss (number of classification errors) of MDL
may be derived with the techniques from
\cite{Hutter:04uaibook}, using the bound on the quadratic
distance. \cite{Hutter:02spupper} gives also bounds for
\emph{arbitrary} loss functions, however this requires a bound
on the Kullback-Leibler divergence rather than the quadratic
distance. Unfortunately, this does not hold for static MDL
\cite{Poland:04mdl}. For the regression setup, analysis of
other, more general or even arbitrary loss functions is even
more demanding and, as far as we know, open.

Considering only discrete model classes is certainly a
restriction, since many models arising in science (e.g. physics or
biology) are continuous. On the other hand there are arguments in
favor of discrete classes. From a computational point of view they
are definitely sufficient. Real computers may even treat only
finite model classes. The class of all programs on a fixed
universal Turing machine is countable. It may be related to
discrete classes of stochastic models by the means of
semimeasures, this is one of the central issues in Algorithmic
Information Theory \cite{Li:97}.

%%%%%%%%%%%%%%%%%%%%%%%%%%%%%%%%%%%%%%%%%%%%%%%%%%%%%%%%%%%%%%%
%         Bibliography        %
%%%%%%%%%%%%%%%%%%%%%%%%%%%%%%%%%%%%%%%%%%%%%%%%%%%%%%%%%%%%%%%

\end{document}